# Reactive Vehicle Guidance using Dynamic Maneuvering Cue

Alexander Von Moll[1] Isaac E. Weintraub[2]

*Air Force Research Laboratory, WPAFB, OH, 45433, USA*

**Recent approaches for navigating among dynamic threat regions (i.e., weapon engagement zones) have focused on planning entire trajectories. Moreover, the allowance for penetration into these threat regions was based on heuristic measurements of risk. This paper offers an approach for a more reactive (i.e., feedback-based) guidance that is based on closed-form analytical expressions and thereby suitable for onboard, real-time execution. In addition, a risk measurement is formulated based upon the concept of Dynamic Maneuvering Cue (DMC), which measures the amount of turn a vehicle would need to take in its current state in order to put itself outside the threat region. This approach is then extended to handle multiple threat regions simultaneously (with minimal additional computational complexity). Finally, the DMC constraint is applied to a simple feedback controller as well as a model predictive controller (MPC). The MPC shows better performance, but at the cost of having to solve an optimization problem online, versus the meager computational burden associated with the simple controller. This approach, which is based on assuming the threats are adversarial, may be used as a conservative method for collision avoidance and deconfliction.**

## I. Introduction

*Urban Air Mobility* (UAM) and *Advanced Air Mobility* (AAM) represent a growing vision for the future of transportation at scale [1], [2]. This concept envisions numerous manned and unmanned systems sharing the same airspace, navigating efficiently to their destinations while remaining free from collision. A central challenge for realizing UAM/AAM is ensuring effective deconfliction and managing collision risk in potentially congested environments. The primary contribution of this work is a method for quantifying the instantaneous risk of collision and the corresponding evasive maneuver required, represented as a Dynamic Maneuvering Cue (DMC), intended to support safe UAM/AAM operations.

### A. The Need for Separation Assurance and Maneuver Guidance

Effective Air Traffic Management (ATM) is fundamental to ensuring the safe and efficient flow of any aircraft, traditional or novel. A core component is separation assurance – preventing aircraft from violating minimum safe distances. Understanding the minimum maneuver required to regain or maintain safe separation is crucial for designing conflict detection and resolution systems [3]. Estimating the probability of collision based on trajectory predictions and uncertainties can inform the urgency and magnitude of necessary avoidance maneuvers [4]. Existing frameworks like the Airborne Collision Avoidance System (ACAS), including its various iterations (TCAS I/II, ACAS Xa/X) [5], [6], aim to reduce mid-air collision probability in traditional aviation. While systems like ACAS X leverage probabilistic models and optimized advisory logic [6], the underlying need remains for clear, actionable guidance, especially as airspace density increases with UAM/AAM.

### B. Fighter Combat Tactics and the Dynamic Maneuvering Cue (DMC) Concept

This paper makes use of concepts from the fighter-pilot community and demonstrates how automation can be performed to improve inter-vehicle traffic management, considering a conservative approach wherein other aircraft are modeled as adversaries. Moreover, the conceptual origins of dynamically

This paper is based on work performed at the Air Force Research Laboratory (AFRL) *Control Science Center* and is supported by AFOSR LRIR #24RQCOR002 (funded by Dr. Frederick Leve). DISTRIBUTION STATEMENT A. Approved for public release: distribution is unlimited (AFRL-2025-5480) Cleared 5 Dec 2025 assessing required maneuvers for survival lie partly in air combat. Fighter pilots must constantly evaluate threat positions and kinematics to execute timely evasive or offensive maneuvers, often based on geometry and relative heading changes needed to defeat a weapon or gain an advantage. The underlying principle of calculating a necessary instantaneous maneuver is central to tactical decision-making. The interested reader can look to a survey of machine learning methods for air-combat behaviors [7].

[1]Aerospace Engineer, Control Science Center.
[2]Senior Electronics Engineer, Control Science Center, AIAA Associate Fellow.

In air combat, one of the many priorities is the avoidance of incoming threats [8]. In an adversarial scenario, the incoming threat aims to capture the target, and the vehicles have limited performance. Aircraft make evasive maneuvers to avoid vehicle collision or capture, the amount of turn that is taken is a measure for how evasive a maneuver must be to avoid capture. The greater the required turn by the aircraft to avoid capture, the greater the risk accumulated by the target aircraft [9]. This concept is known by pilots as a *Dynamic Maneuver Cue* (DMC). In this work the DMC is defined as an indication of the number of degrees of heading change required to evade capture by an adversarial pursuer.

The greater the DMC, the greater the maneuver required to evade capture. DMC can be used as a measure of risk by the target aircraft, as a DMC of 0 degrees simply states that no heading change is required to evade capture – the target aircraft is at no risk from the adversarial pursuer. If the DMC were 180 degrees, then an about hard turn is required to evade capture and the risk taken by the platform is quite significant. A DMC of 180 degrees means that the aircraft is in danger of capture no matter the maneuver it takes. The establishment of rules that pilots can use to make timely decisions is used in practice and was even investigated in a work by Burgin and Sidor [10].

### C. Related Approaches

The idea of calculating required maneuvers based on geometry is not unique to fighter tactics. Research in vehicle deconfliction and collision avoidance often focuses on algorithms that compute necessary heading or velocity changes. For instance, the Velocity Obstacle (VO) method provides a geometric approach to determine the range of velocities (including headings) that would lead to a collision, thereby informing the selection of a safe velocity [11]. The leveraging of control barrier functions has shown promise for addressing a scalable way to ensure safe separation distance between aircraft with turn-constraints [12]. Other work focuses on coordinated conflict resolution, such as automated methods for air traffic control in terminal areas aimed at minimizing delays while ensuring separation [13]. Planning paths to avoid keep-out-zones or obstacles is a broad field, some techniques such as vector fields [14], potential fields [15], sampling-based methods [16], and specialized approaches like terrain avoidance [17] are particularly related to this work. Developing effective sense-and-avoid (SAA) capabilities, potentially using onboard sensors like radar [18], is critical for UAS integration [19].

Furthermore, the concept of navigating safely relative to potential threats relates to path planning around Engagement Zones (EZ), particularly the Basic Engagement Zone (BEZ) [20]. The BEZ framework considers an opponent's "worst-case" strategy, assuming they will maneuver optimally to intercept. This conservative, adversarial perspective is analogous to ensuring separation from other airspace users who may not be cooperative or predictable. This paper specifically leverages the geometric formulation of the BEZ [20] as the foundation for calculating the DMC, treating other vehicles as potential worst-case threats to ensure a high degree of safety.

### D. Risk Management and Heuristic Control

Quantifying the required maneuver via the DMC directly informs risk management. Concepts like Collision Risk Modeling (CRM) assess accident probability based on factors like traffic density and separation standards [21], [22]. The DMC provides a real-time, state-dependent measure related to this risk – a larger DMC indicates a situation demanding a more significant deviation, reflecting higher immediate risk. While optimal control solutions for trajectory planning exist, real-time safety often relies on heuristic controls and decision aids that provide rapid, understandable guidance based on simplified models or expert knowledge [23]. A clearly defined DMC, as developed in this work, is well-suited to serve as a key input for such systems, translating complex state information (relative positions, velocities) into a simple, actionable cue for pilots or autopilots to initiate timely avoidance actions.

### E. Outline

The main contributions in this work are as follows:

- A mathematical description of Dynamic Maneuvering Cue
- An approach for considering multiple adversarial systems using DMC
- Control strategies based on DMC.

This work begins with Technical Preliminaries in Section II. Simulation and results of some example scenarios are in Section V. Final conclusions and future work appears in Section VI.

## II. Technical Preliminaries

In this section, we define the models used for the vehicle as well as the model for the engagement zone associated with the other vehicles. Let $A = (x_A, y_A)^\top \in \mathbb{R}^2$ denote the position of the vehicle being controlled (i.e., the Agent) in the 2D plane. The Agent moves with *simple motion*, that is, its kinematics are given by

$$\dot{A} = (\dot{x}_A, \dot{y}_A)^\top = (v_A \cos \psi, v_A \sin \psi)^\top, \qquad (1)$$

where $v_A$ is the constant speed of the Agent and $\psi$ is its heading angle w.r.t. the positive $x$-axis. Now consider another vehicle which potentially poses a risk of collision with $A$; let this vehicle be denoted $T$ (for Threat). $T$ has a speed $v_T > v_A$ for which we define a speed ratio $\mu = \frac{v_A}{v_T} < 1$. Since $A$ wishes to avoid collision with $T$, a conservative approach would be to assume that $T$ will put itself on a collision course with $A$ (given $A$'s current course, $\psi$) over the next period of time. Denote this time as $t_r$ – it represents a reaction time of some sort which may be related to the actual vehicle dynamics and performance of $A$ (which are not modeled explicitly here). The reaction time specifies a maximum distance that $T$ can travel of

$$R = v_T t_r. \qquad (2)$$

Next consider a minimum safety distance $r$ which represents the minimum distance $A$ can be from $T$ without colliding (possibly plus an additional safety margin).

Based on these parameters, the *basic engagement zone* (BEZ) of $T$ w.r.t. $A$ is the set of positions for $A$ such that it is possible for $T$ to collide with $A$ within $t_r$ time if the latter holds course. This model is adapted from [20] in which $T$ represented a weapon with a maximum range of $R$. Note that $T$'s heading doesn't factor in to the definition of the BEZ as it is implicitly assumed that it, too, moves with simple motion (i.e., can turn instantaneously). From [20], we have that the BEZ distance (from the threat's position), as a function of aspect angle, is given by

$$\rho(\xi) = \mu R \left[ \cos \xi + \sqrt{\cos^2 \xi - 1 + \frac{(R+r)^2}{\mu^2 R^2}} \right]. \qquad (3)$$

Let the distance from $A$ to $T$ be denoted by $d$. When $d \leq \rho$, $A$ is inside $T$'s BEZ. As shown in [20], for a particular heading $\psi$, the locus of points that define the boundary of the BEZ (i.e., where $d = \rho$) is given by a circle of radius $R + r$ whose origin, $B$, is offset a distance $\mu R$ from $T$'s position in the direction of $-\psi$. Fig. 1 shows an example BEZ (green). In the Fig. 1, the agent, $A$, moves with constant speed, $v_A$, in the direction of the blue arrow. The solid red circle indicates $T$'s capture radius, $r$, the smaller dotted red circle is $T$'s reachability region, and the larger is $T$'s capturability region (i.e., all the positions that can be contacted by $T$'s capture circle).

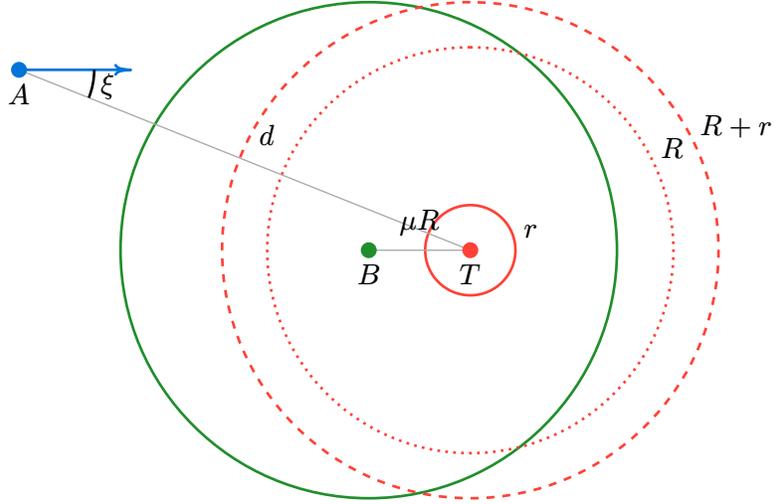

**Fig. 1  Basic engagement zone (adapted from [20])**

### III. Dynamic Maneuvering Cue

#### A. Minimum DMC to get outside the BEZ

We wish to find the aspect angle, $\xi^*$, for a vehicle inside of a BEZ such that the vehicle *would* be on the boundary of the BEZ associated with its new heading. Then, the difference between the vehicle's current heading and this particular aspect angle would correspond to the Dynamic Maneuvering Cue (DMC) value. That is, how much the vehicle would need to turn in order to no longer be inside the threat's BEZ.

Now let $d$ be defined as the distance between the threat and vehicle positions. The 'minimally safe' aspect angle, $\xi^*$ corresponds to the one in which $\rho(\xi^*) = d$:

$$d = \mu R \left[ \cos \xi + \sqrt{\cos^2 \xi - 1 + \frac{(R+r)^2}{\mu^2 R^2}} \right]$$

$$d - \mu R \cos \xi = \mu R \sqrt{\cos^2 \xi - 1 + \frac{(R+r)^2}{\mu^2 R^2}} \tag{4}$$

$$d^2 - 2\mu R d \cos \xi + \mu^2 R^2 \cos^2 \xi = \mu^2 R^2 \left( \cos^2 \xi - 1 + \frac{(R+r)^2}{\mu^2 R^2} \right)$$

$$d^2 - 2\mu R d \cos \xi = -\mu^2 R^2 + (R+r)^2$$

which yields

$$\xi^* = \cos^{-1} \left( \frac{d^2 + \mu^2 R^2 - (R+r)^2}{2\mu R d} \right). \tag{5}$$

Note that there are always two solutions for the above expression. Let us denote these two solutions as $\xi_1^*$ and $\xi_2^*$, and, by definition, $\xi_2^* = -\xi_1^*$. Therefore, when $\xi \in \xi_{\text{unsafe}} \equiv [-\xi^*, \xi^*]$, $A$ is inside the BEZ. Also, $\xi^*$ is only defined when the following condition is satisfied

$$\left| d^2 + \mu^2 R^2 - (R+r)^2 \right| \leq 2\mu R d \tag{6}$$

since the argument to the inverse cosine must be in $[-1, 1]$. When this condition is not satisfied, we define $\xi_{\text{unsafe}}$ as the empty set. Finally, the DMC value is given by,

$$\text{DMC}(\xi, d; \mu, R, r) = \begin{cases} \text{sign}(\xi) \min_{\xi_c \in \{-\xi^*, \xi^*\}} |\xi - \xi_c|, & \text{if } \xi \in \xi_{\text{unsafe}} \\ 0, & \text{otherwise.} \end{cases} \tag{7}$$

where $\xi$ is the current aspect angle of the vehicle w.r.t. the line of sight from the threat.

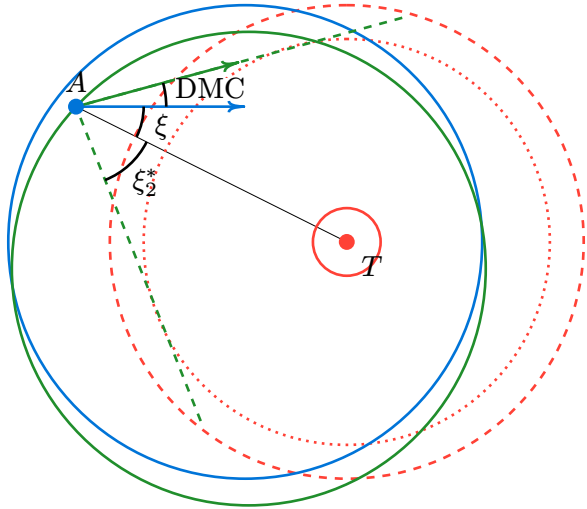

**Fig. 2  Pictorial representation of the analytic DMC. $\mu = 0.5, R = 0.9, r = 0.15$**

Fig. 2 shows the results of this process for a particular set of parameters. Note that the resulting heading puts $A$ exactly on the boundary of the BEZ, but its heading may still point *into* the BEZ.

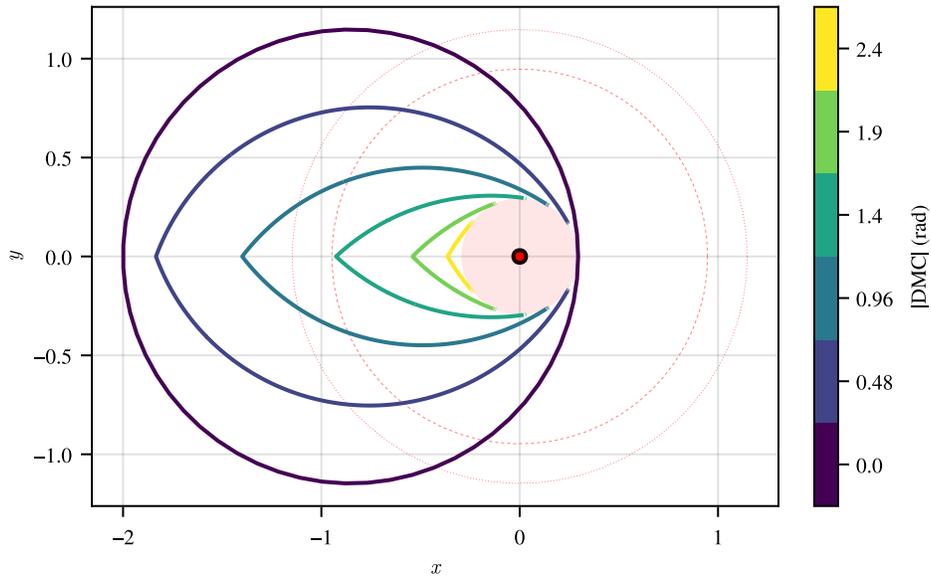

**Fig. 3  DMC absolute values (rad) for $A$ heading to the right ($\psi = 0$).**

Fig. 3 shows the absolute values of DMC as a function of $A$'s position for a heading of $\psi = 0$. Note that the DMC = 0 contour corresponds to the boundary of the BEZ and the DMC outside of the BEZ is also zero. The shaded circle centered on the origin is a circle whose radius is

$$d_{\text{no escape}} = (1 - \mu)R + r, \qquad (8)$$

which corresponds to positions from which $A$ cannot guarantee escape from $T$ under any heading (c.f., [24] for more details).

### B. DMC for staying outside BEZ

In order to prevent $A$ from going into the BEZ immediately, one may consider a heading for which the point $A$ is *tangent* to the BEZ. As a result of the requirement for tangency, a right triangle arises. The 'minimally safe' aspect angle, in this case, is

$$\xi^* = \sin^{-1}\left(\frac{R+r}{d}\right). \tag{9}$$

As before, there are two solutions. When $A$ is sufficiently far from $T$, the BEZ model implies that if $A$ heads towards the tangent of the capturability region that $A$ also heads towards the tangent of the BEZ as shown in Fig. 4.

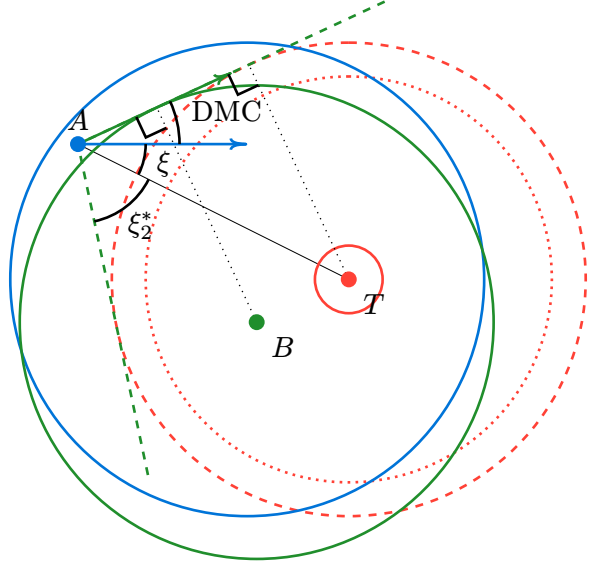

**Fig. 4  DMC for staying outside BEZ - the tangency DMC.** $\mu = 0.5, R = 0.9, r = 0.15$.

Based on Eq. (9), this will only work when $d \geq R + r$ since the argument to inverse sin must be in the range $[-1, 1]$. This corresponds to positions in which $A$ is outside of the capturability region.

As stated previously, when $A$ is sufficiently far from $T$ these two conditions are coincident: the minimally safe aspect angle is the one in which $A$ is tangent to both the capturability region and the BEZ. As $A$'s position is moved closer to $T$, eventually these two conditions cannot be satisfied simultaneously. Fig. 5 shows one such critical case.

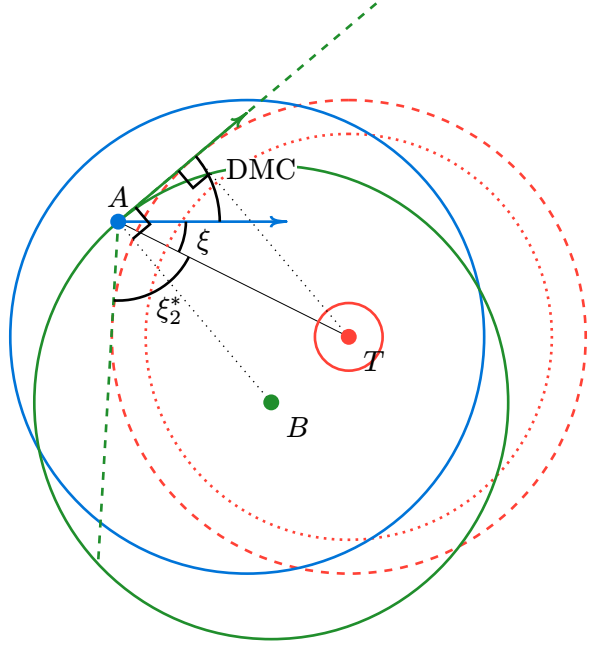

**Fig. 5 Critical case where $A$ is simultaneously tangent to the BEZ and to the capturability region. $\mu = 0.5, R = 0.9, r = 0.15$**

The critical distance occurs when $\triangle TAB$ is a right triangle. Thus,
$$d_{\text{crit}} = \sqrt{(R+r)^2 + \mu^2 R^2}. \tag{10}$$
Below the critical distance, it is not possible to choose a heading to be on the boundary of the BEZ *and* tangent to it simultaneously. Therefore, the capturability tangency condition is sufficient for staying outside the BEZ but not necessary.

The following expressions summarize the prescribed DMC computation for the case of staying outside the BEZ (notated with overlines to differentiate them from the previously defined quantities).

$$\bar{\xi}^* = \begin{cases} \cos^{-1}\left(\frac{d^2 + \mu^2 R^2 - (R+r)^2}{2\mu R d}\right), & \text{if } d < d_{\text{crit}} \\ \sin^{-1}\left(\frac{R+r}{d}\right), & \text{otherwise.} \end{cases} \tag{11}$$

$$\overline{\text{DMC}}(\xi, d; \mu, R, r) = \begin{cases} 0, & \text{if } d > \rho(\xi; u, R, r) \\ \text{sign}(\xi) \min_{\xi_c \in \{-\xi^*, \xi^*\}} |\xi - \xi_c|, & \text{otherwise.} \end{cases} \tag{12}$$

### C. Extension to Multiple Threats

The concept of DMC based on the BEZ can be readily extended to the case of multiple Threats. Now consider Threats $T_i, i = 1, ..., N$. Each Threat $T_i$ gives rise to a particular safe cone of headings, which, if implemented by $A$, would place $A$ outside of $T_i$'s BEZ. This safe cone is defined by the following
$$\Psi_i = [\lambda_i + \xi^*, \lambda_i - \xi^*] \tag{13}$$
where $\xi^*$ is defined either by Eq. (5) or Eq. (11) and the line of sight angle is defined by
$$\lambda_i = \operatorname{atan}\left(y_{T_i} - y_A, x_{T_i} - x_A\right). \tag{14}$$
The cone of headings for $A$ which are safe w.r.t. *all* of the Threats is given by the set intersection of individual cones:
$$\Psi = \bigcap_i^N \Psi_i = \left[\underline{\psi}, \overline{\psi}\right] \tag{15}$$

Note that special care must be taken with the ranges, $\Psi_i$, and set operations due to the wrapping of the angular quantities. Finally, the overall DMC is found by finding the minimum angular distance from $A$'s current heading, $\psi$, to one of the boundaries of $\Psi$ (if $\psi$ is not already safe):

$$\psi_c = \underset{\hat{\psi} \in \{\underline{\psi}, \overline{\psi}\}}{\operatorname{argmin}} \left| \psi - \hat{\psi} \right|$$

$$\mathrm{DMC}(\psi; \mu, R, r) = \begin{cases} 0, & \text{if } \psi \in \Psi \\ \psi - \psi_c, & \text{otherwise.} \end{cases} \quad (16)$$

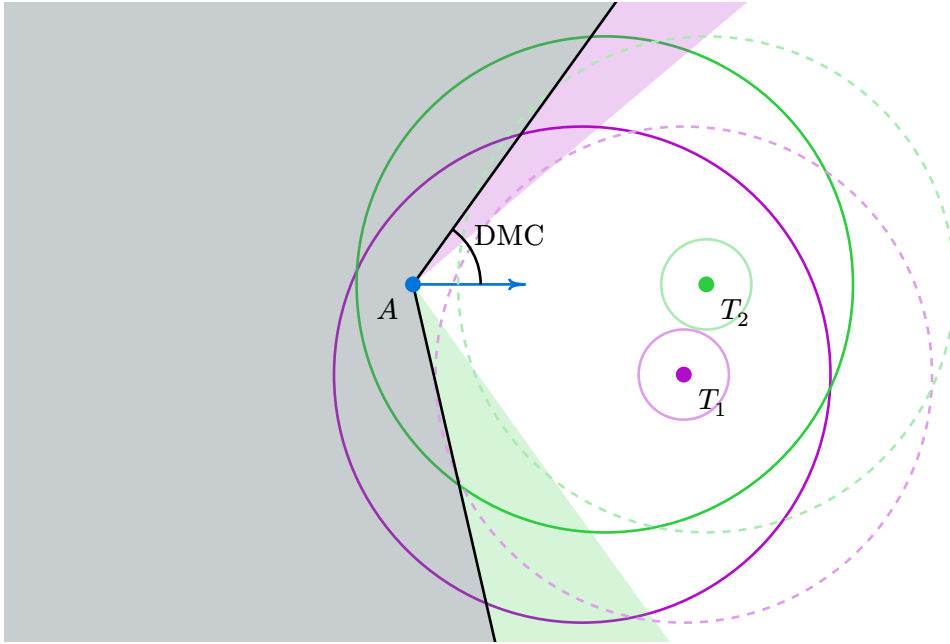

**Fig. 6  Pictorial representation of safe cone of headings (i.e., the overlap region delineated with black lines) and DMC for multiple Threats.**

## IV. Controllers

In this section, we apply the concept of DMC to a navigation example. Define a goal point $G \in \mathbb{R}^2$ which the Agent wishes to reach all while respecting a constraint on the DMC with respect to one or more Threats. In general, $A$ may not necessarily be required to maintain DMC = 0 along the entire trajectory. Rather, $A$ is assumed to have a DMC threshold denoted by $\varepsilon$, i.e., the following path constraint is applied for all the subsequently defined controllers:

$$g(\mathbf{x}(t)) = |\mathrm{DMC}| - \varepsilon \leq 0, \quad \forall t \in [0, t_f] \quad (17)$$

where $\mathbf{x}$ is the full state of the system (i.e., the stacked vector of positions of the vehicles as well as $A$'s heading) and $t_f$ is the first time instant in which $A$ and $G$ are coincident.

### A. Simple (Aim at Goal)

Perhaps the most simple control scheme for the navigation example (other than setting a constant heading) is to select a heading at each time instant in the direction towards the goal, $G$. That is, the nominal control is

$$\psi_{\mathrm{nom}}(A, G) = \operatorname{atan}(y_G - y_A, x_G - x_A), \quad (18)$$

which is a feedback control law based on the instantaneous position of $A$. When the DMC constraint, Eq. (17), becomes active, then $A$ can no longer apply $\psi_{\mathrm{nom}}$ and must instead choose a heading which respects the constraint. In these cases, without planning into the future at all to determine which heading to take, it is reasonable to find the nearest heading to the desired $\psi_{\mathrm{nom}}$ which still satisfies the DMC threshold. The controller can be written as

$$\psi_s = \begin{cases} \psi_{\text{nom}}, & \text{if } |\text{DMC}| \leq \varepsilon \\ \psi_{\text{nom}} + \text{DMC} - \text{sign}(\text{DMC})\varepsilon, & \text{otherwise,} \end{cases} \quad (19)$$

where the subscript $s$ denotes "simple".

One potential pitfall of the Simple controller is that there are configurations for which the nominal desired heading, $\psi_{\text{nom}}$ lies directly in the center of the unsafe cone of headings. We refer to such a configuration as lying on a switching surface (SS) which is defined as

$$\mathcal{S}(T_1, ..., T_N, G) = \left\{ A \, \big| \, |\psi_{\text{nom}} - \underline{\psi}| = |\psi_{\text{nom}} - \overline{\psi}| \right\}, \quad (20)$$

i.e., the desired heading is equidistant to either boundary of the safe cone of headings. This essentially is the dividing line between heading CCW or CW around one or more Threats.

In the case of a single Threat, $\mathcal{S}$ is actually a dispersal surface meaning that if the state of the system started on $\mathcal{S}$ it would immediately depart it in forward time, never to return. However, in the case of two or more Threats, $\mathcal{S}$ (or at least parts of it) are "sticky". That is, if the state of the system lies on $\mathcal{S}$, the Simple controller will push the state of the system onto the opposite side of $\mathcal{S}$ wherein the Simple controller will immediately push the state back across $\mathcal{S}$. Generally, each Threat contributes a non-sticky branch to $\mathcal{S}$ while each pair of Threats can produce a sticky branch.

### B. Model Predictive Control

One option for increasing the sophistication of the control approach is to use model predictive control (MPC) [25] which, in some sense, bridges the gap between path planning and control. Instead of reacting strictly to the current state, as in the case of the Simple controller, an MPC-based controller optimizes a sequence of control inputs over a finite time horizon and then implements some subsequence thereof (often just the first control input in the optimized sequence). What follows is a formulation of a very simple MPC based on DMC.

It is assumed here that the Threat location(s) are stationary, however, this formulation could be modified to account for moving Threat(s) if their motion was known or could be estimated. In lieu of modeling the motion of the Threat(s) explicitly, one can simply recompute either the Simple or MPC controller at every time instant and use the updated Threat position(s). Let $t_s$ be the sample time and $t_h = H t_s$ be the time horizon for planning, where $H+1$ is the number of samples in the planning horizon. Denote the sequence of control inputs over the planning horizon by $\boldsymbol{\psi} = (\psi(0), \psi(t_s), \psi(2t_s), ..., \psi(t_h))^\top$. The following optimization problem encourages the Agent to end up closer to its goal, $G$, while respecting its prespecified DMC threshold, $\varepsilon$:

$$\boldsymbol{\psi}^* = \underset{\boldsymbol{\psi}}{\text{argmin}} \; \overline{A(t_h)G} \quad (21)$$
$$\text{subject to } |\text{DMC}|\big|_{t=kt_s} - \varepsilon \leq 0, \quad k = 0, ..., H.$$

As will be shown in the results to follow, the MPC approach is less myopic than the Simple controller and can often outperform it as well as avoid getting trapped in troublesome regions of the state space where the latter can get stuck. This increased robustness and performance comes at the cost of increased computation since solving Eq. (21) requires running a nonlinear program solver (NLP) versus the closed form expression in Eq. (19) for the Simple controller.

### V. Simulations & Results

In this section, the controllers described in the preceding section are simulated in a variety of scenarios to demonstrate how they operate and the types of scenarios in which they are applicable. For each of the simulations, it is assumed that the Threat does not 'fire' and/or begin pursuing the Agent unless otherwise noted. Therefore, the resulting trajectories are optimistic in the sense that $A$ is not forced to make an evasive maneuver.

### A. Comparison to Optimal Control

Prior work [20] focused on development of the BEZ model and included a path planning example around a single Threat. There, the problem was formulated as an optimal control problem, i.e., to reach $G$ in minimum time without entering inside the BEZ. This corresponds to setting the DMC threshold to $\varepsilon = 0$.

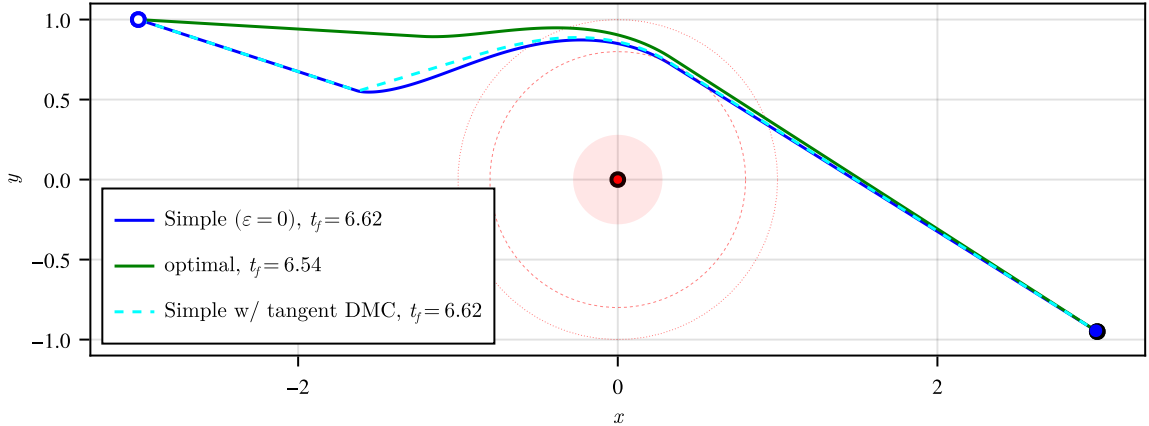

**Fig. 7   Comparison to optimal from [20].**

Fig. 7 shows an example where $A$ must go around $T$ (which is located at the origin) in order to reach its goal, $G$. For this example the parameters are set to $v_A = 1$, $\mu = 0.9$, $r = 0.2$, and $R = 1 - r$. The Simple controller was successful in guiding the Agent to its goal while respecting the DMC constraint. Moreover, it took only 1.13% more time to reach the goal compared to the optimal control. Usage of the tangent DMC, Eq. (11), produces a similar trajectory, as shown in Fig. 7, but exhibits an abrupt heading change at the first instant that the Agent reaches the BEZ.

### B.  Multiple Threats and Varying Thresholds

One of the major advantages of the DMC concept is that it provides a physics-based approach for quantifying the risk for a particular vehicle against a potential threat. One prior work [26] allowed for penetration into the engagement zone of a Threat by formulating an objective cost functional which mixed time-to-go with a notion of risk accrued over the trajectory. The risk model used there had the form of

$$\text{risk}(\mathbf{x}(t), \psi(t)) = \begin{cases} \frac{\rho}{d} - 1, & d \leq \rho \\ 0, & \text{otherwise}, \end{cases} \quad (22)$$

where $\rho$ is the distance of the Threat location to the boundary of its engagement zone for the aspect angle arising from $\psi$ (e.g., Eq. (3)). That model's value has no inherent meaning – it simply starts at 0 on the boundary of the engagement zone and increases monotonically as the distance to the Threat decreases. In contrast, the DMC value directly measures the magnitude of the evasive maneuver required for $A$ to guarantee safety.

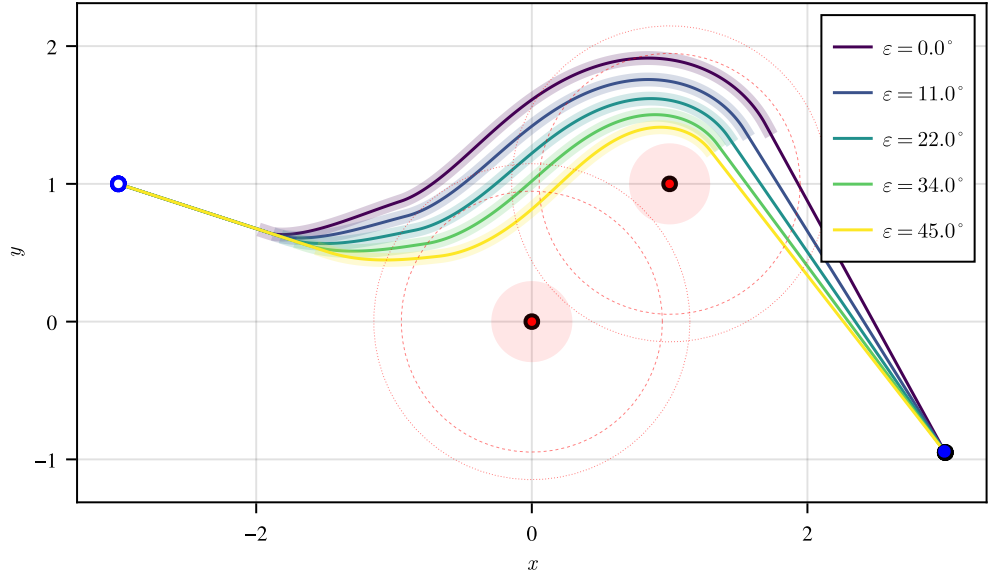

**Fig. 8  Simple controller applied to scenario with 2 Threats with various DMC thresholds.**

Fig. 8 shows several trajectories simulated using the Simple controller against 2 Threats with the DMC threshold, $\varepsilon$, set to different values (note: $R = 0.947$). The trajectories with higher DMC thresholds navigate closer to the Threats thereby reaching the goal location in shorter time. With a non-zero DMC threshold, $A$ is able to navigate inside of the BEZs (c.f., Fig. 3). The highlighted portions of the trajectories in Fig. 8 indicate where the DMC constraint is active. Note that going clockwise around the Threats is clearly suboptimal in this case which further highlights the myopic nature of the Simple controller.

### C. Model Predictive Control

The MPC-based controller is demonstrated on the same example as in the previous section. Note that the initial conditions are such that $A$ starts slightly above the line of symmetry w.r.t. the closer BEZ and the goal. That is why the Simple controller results in heading up once the DMC becomes nonzero. In contrast, the MPC-based controller (with $H = 25$ and sample time $t_s = 0.07$) is able to predict far enough ahead to see that going upward results in less progress towards the goal because of the second BEZ (see Fig. 9).

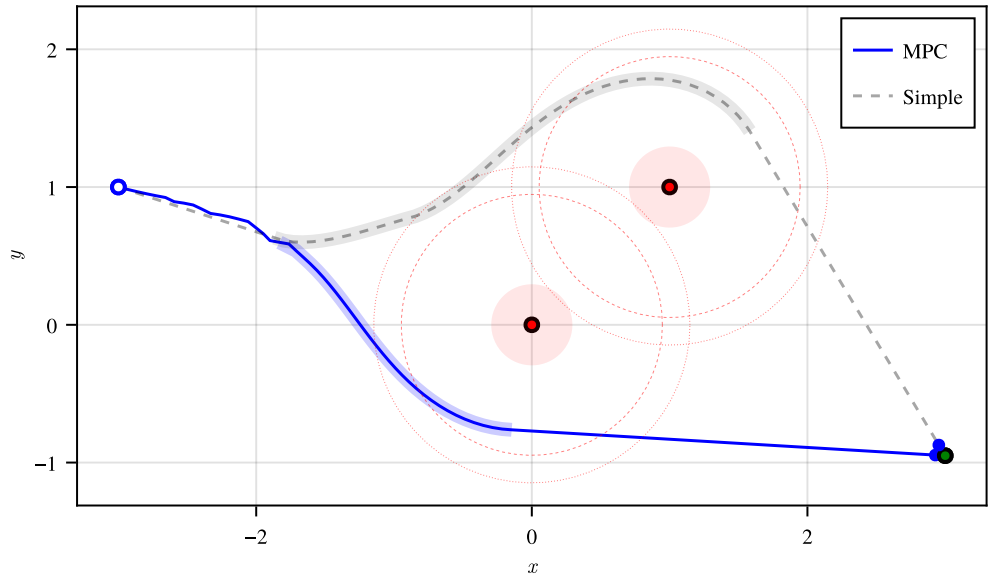

**Fig. 9  MPC-based controller results in shorter path around multiple BEZs ($\varepsilon = 10°$).**

Of course, this comes with the additional computational burden of solving the optimization in Eq. (21) at each timestep. There are also some numerical challenges with the possibility for the DMC to change discontinuously when multiple BEZs are considered. It is possible, for example, for the joint safe range of headings, Eq. (15), to become noncontiguous. Such challenges partly explain the nonsmooth portion of the MPC trajectory.

### D. Moving Threat

This last result is a demonstration of the straightforward application of the Simple controller (and DMC, in general) to a scenario in which the threat, itself, is moving. In this case, $T$ is assumed to be a vehicle carrying a weapon onboard which gives rise to the associated BEZ. Instead of the BEZ origin remaining stationary, it remains fixed to the adversarial vehicle's position. Fig. 10 shows an example wherein $T$ employs pure pursuit and moves with a speed $v_T = 0.6$ (compared to the Agent's speed of $v_A = 1$).

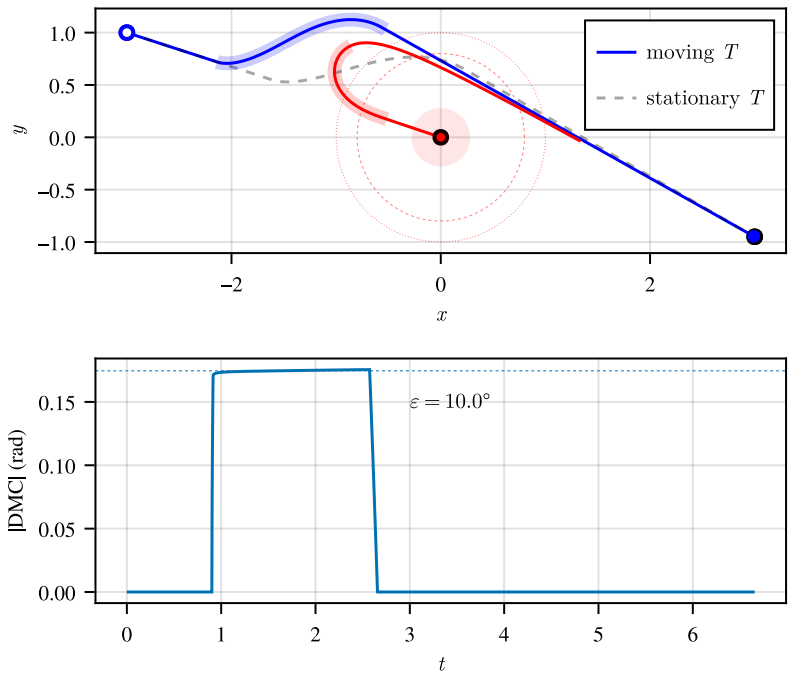

Fig. 10  Simple controller against a mobile threat moving with $v_T = 0.6$ using pure pursuit.

Because $T$ actively pursues $A$, the latter must begin its avoidance maneuver much earlier compared to the case wherein $T$ is stationary. This is simply because the DMC constraint becomes active earlier as well – $A$'s controller makes no consideration of $T$'s motion or strategy, it simply reacts to the instantaneous threat posture.

## VI. Conclusions

Risk quantification is an important aspect of vehicle guidance. This paper presented a novel physics-based risk metric for vehicles in the presence of threats. The Dynamic Maneuvering Cue (DMC) indicates the amount of turn the vehicle must make in order to avoid the possibility for future collision/intercept by the threat. Higher DMC values indicate the need for a more aggressive maneuver which translates into higher risk. Although the DMC value is associated with a worst-case maneuver by the threat, it may still serve as a conservative risk measure in contexts where the 'threats' are not necessarily adversarial. Two different control approaches were developed based on treating the DMC as a constraint: a simple feedback controller which aims as close to the goal location as possible, and an MPC-based controller. The MPC-based controller avoids some of the pitfalls of the simple controller but at a much higher computational cost.

The DMC introduced in this paper assumes that the vehicle moves with a constant speed and does not take into consideration its performance characteristics. Future work will focus on extending the DMC concept to incorporate more information about the vehicle and its dynamics. For example, in order to avoid the threat region, the vehicle may have to initiate a maximum turn-rate turn and/or speed up or slow down. Additionally, direct comparison and/or incorporation to other collision avoidance approaches (including control barrier functions) is left for future work.